\documentclass[10pt, letterpaper,reqno]{amsart}
\usepackage{amssymb}
\usepackage{amsmath, dsfont}
\usepackage{amsfonts, color}
\usepackage{graphicx}
\usepackage{mathrsfs,amssymb, slashed, cite}

\usepackage{hyperref}
\usepackage{cleveref}

\newcommand{\R}{\mathbb{R}}

\renewcommand{\H}{\mathcal{H}}
\newcommand{\1}{\mathds{1}}

\newtheorem{definition}{Definition}

\newtheorem{proposition}{Proposition}[section]

\newtheorem{theorem}{Theorem}[section]
\numberwithin{equation}{section}

\newcommand{\qtq}[1]{\quad\text{#1}\quad}

\allowdisplaybreaks

\begin{document}

\title[Endpoint Strichartz]{The double-endpoint Strichartz estimate for dispersion-managed Schr\"odinger equations}

\author[J. Murphy]{Jason Murphy}
\address{Department of Mathematics, University of Oregon, Eugene, OR, USA}
\email{jamu@uoregon.edu}

\author[J. Zheng]{Jiqiang Zheng}
\address{Institute of Applied Physics and Computational Mathematics and National Key Laboratory of Computational Physics, Beijing 100088, China}
\email{zheng\_jiqiang@iapcm.ac.cn}
\thanks{J. Zheng is supported
by National Key R\&D program of China (2021YFA1002500) and NSF Grant of China (No. 12271051).}

\begin{abstract} We adapt the near-far argument of Hadama and Toyoshima to establish the double-endpoint Strichartz estimate for dispersion-managed Schr\"odinger equations.
\end{abstract}

\maketitle

The recent paper \cite{Hadama} provides a new proof of the Keel--Tao double-endpoint Strichartz estimate \cite{KeelTao, Foschi}, based on the near-far argument introduced in \cite{HadamaToyoshima}.  In this note we show that this argument can also be adapted to treat dispersion-managed Schr\"odinger equations\footnote{See e.g. \cite{Hundertmark, Turitsyn, MVH} for physical and mathematical introductions to dispersion-managed Schr\"odinger equations.} of the form
\begin{equation}\label{LS}
i\partial_t u + \gamma(t)\Delta u = 0.
\end{equation}

For a suitable class of dispersion maps $t\mapsto\gamma(t)$, the work \cite{MVH} established Strichartz estimates for \eqref{LS}, but left the double-endpoint estimate open, due to the fact the Christ--Kiselev argument breaks down in this case.  While the fact that the total dispersion map $\Gamma(t,s)=\int_s^t\gamma(\tau)\,d\tau$ may vanish away from the diagonal $s=t$ prevents the direct application of existing estimates, we will show that for suitable dispersion maps the arguments of \cite{Hadama} may be suitably adapted to address this missing case.  As in \cite{MVH}, we will consider dispersion maps satisfying the following conditions:
\begin{definition} We call $\gamma:\R\to\R$ \emph{admissible} if it satisfies the following:
\begin{itemize}
\item $\gamma$ is one-periodic: $\gamma(t+1)=\gamma(t)$ for all $t\in\R$.
\item $\gamma$ and $\gamma^{-1}$ are bounded: $\|\gamma\|_{L^\infty}+\|\gamma^{-1}\|_{L^\infty} <\infty$.
\item $\gamma$ has at most finitely many discontinuities on $[0,1]$.
\item $\gamma$ has nonzero average over its period: $\langle \gamma\rangle:=\int_0^1 \gamma(t)\,dt \neq 0.$
\end{itemize}
\end{definition}

We will prove the following:
\begin{theorem}[Double-endpoint Strichartz estimate]\label{T} Fix $d\geq 3$. Let $\gamma$ be an admissible dispersion map and define
\[
\Gamma(t,s)=\int_s^t \gamma(\tau)\,d\tau. 
\]
Then there exists $C=C(d,\gamma)>0$ such that
\begin{equation}\label{de}
\biggl\| \int_{-\infty}^t e^{i\Gamma(t,s)\Delta}F(s)\,ds\biggr\|_{L_t^2 L_x^{\frac{2d}{d-2}}(\R\times\R^d)} \leq C\|F\|_{L_t^2 L_x^{\frac{2d}{d+2}}(\R\times\R^d)}.
\end{equation}
\end{theorem}

As is shown in \cite[Proof of Theorem~1.1]{Hadama}, one can use a duality argument together with a suitable factorization of the input $F$ to derive \eqref{de} from the following bound:
\begin{equation}\label{de1}
\biggl\| \int_{-\infty}^t f(t) e^{i\Gamma(t,s)\Delta}[g(s) F(s)]\,ds\biggr\|_{L_{t,x}^2} \lesssim \|f\|_{L_t^\infty L_x^{d,\infty}}\|g\|_{L_t^\infty L_x^{d,\infty}}\|F\|_{L_{t,x}^2}.
\end{equation} 
Indeed, \eqref{de1} yields a Lorentz-improved version of \eqref{de}.  Thus we focus on establishing \eqref{de1}.

In \cite[Step 4]{Hadama}, the estimate \eqref{de1} is further reduced to the case $g=f\geq 0$, with 
\begin{equation}\label{form}
f^2 = \sum_{k\in S} 2^k \1_{E_k},
\end{equation}
where $E_k\subset\R\times\R^d$ are pairwise disjoint and $\# S<\infty$. Let us utilize this further reduction as well.

Writing 
\[
E_k(t) = \{x:(t,x)\in E_k\}\qtq{and}f(t)=\sum_{k\in S} \1_{E_k(t)} f(t),
\]
the \emph{near-far} decomposition takes the following form: given $F\in L_{t,x}^2$, we write
\begin{align*}
\int_{-\infty}^t f(t)&  e^{i\Gamma(t,s)\Delta} [f(s)F(s)]\,ds  = \mathcal{I}_{\text{near}}+\mathcal{I}_{\text{far}} \\
& := \sum_{k,k'\in S} \int_{\mathcal{R}_{\text{near}}} \1_{E_k(t)}f(t) e^{i\Gamma(t,s)\Delta} [f(s)\1_{E_{k'}(s)}F(s)]\,ds \\
& \quad + \sum_{k,k'\in S}\int_{\mathcal{R}_{\text{far}}} \1_{E_k(t)}f(t) e^{i\Gamma(t,s)\Delta}[ f(s)\1_{E_{k'}(s)}F(s)]\,ds,\end{align*}
where we define the near and far regions by
\[
\mathcal{R}_{\text{near}}=\mathcal{R}_{\text{near}}(t;k,k') = \{s\in(-\infty,t): |\Gamma(t,s)|\leq \ell(k,k')\}
\]
and
\[
\mathcal{R}_{\text{far}}=\mathcal{R}_{\text{far}}(t;k,k')=\{s\in(-\infty,t):|\Gamma(t,s)|>\ell(k,k')\},
\]
and $\ell(k,k')$ will be determined below.  We note that in the case $\gamma(t)\equiv 1$ (so that $\Gamma(t,s)=t-s$) this decomposition reduces to the one used in \cite{Hadama}.

To keep the formulas within the margins, it will be useful to introduce the notation
\[
H_{k,k'}(t,s) = \1_{E_k(t)}f(t) e^{i\Gamma(t,s)\Delta}[f(s)\1_{E_k'(s)}F(s)].
\]

We begin with the observation that for $*\in\{\text{near},\text{far}\}$,
\begin{equation}\label{l1l2}
\|\mathcal{I}_{*}\|_{L_{t,x}^2} \leq \biggl\| \int_{\mathcal{R}_*} H_{k,k'}(t,s)\,ds\biggr\|_{\ell_k^2 \ell_{k'}^1 L_{t,x}^2}. 
\end{equation}
This would follow immediately from the triangle inequality if the right-hand were estimated in the $\ell_{k,k'}^1 L_{t,x}^2$-norm. To obtain the $\ell_k^2$-norm, we rely on the disjointness of the sets $E_k$. In particular, writing $\mathcal{H}_{k,k'}$ for $\int_{\mathcal{R_*}} H_{k,k'}\,ds$, 
\begin{align*}
\biggl\| \sum_{k,k'} \1_{E_k} \H_{k,k'} \biggl\|_{L_{t,x}^2}^2 &= \iint \biggl| \sum_k \1_{E_k} \biggl[\sum_{k'} \H_{k,k'}\biggr] \biggr|^2\,dx\,dt \\
& = \iint \sum_k \1_{E_k} \biggl| \sum_{k'} \H_{k,k'}\biggr|^2\,dx\,dt \quad\text{(using disjointness)}\\
& = \sum_k \biggl\| \1_{E_k} \sum_{k'} \H_{k,k'}\biggr\|_{L_{t,x}^2}^2 \leq \| \H_{k,k'}\|_{\ell_k^2 \ell_{k'}^1 L_{t,x}^2}^2.
\end{align*}

We first consider the near region.  We will use the fact that
\[
\sup_{t_1,t_2\in\R} \|e^{i\Gamma(t_1,t_2)\Delta}\|_{L_x^2\to L_x^2} = 1.
\]

For $k\geq k'$, we use the triangle inequality and H\"older's inequality to estimate 
\begin{align*}
\biggl\| \int_{\mathcal{R}_{\text{near}}} H_{k,k'}(t,s)\,ds\biggr\|_{L_{t,x}^2} & \leq \bigl\| \1_{\mathcal{R}_{\text{near}}(t;k,k')}(s) \|H_{k,k'}(t,s)\|_{L_x^2} \bigr\|_{L_t^2 L_s^1} \\
& \leq \|\1_{\mathcal{R}_{\text{near}}(t;k,k')}(s)\|_{L_t^\infty L_s^2} \|H_{k,k'}(t,s)\|_{L_{t,s,x}^2}.
\end{align*}
Now, writing $e^{i\Gamma(t,s)\Delta}=e^{i\Gamma(t,0)\Delta}e^{i\Gamma(0,s)\Delta}$, we can estimate
\begin{align*}
\|H_{k,k'}(t,s)\|_{L_{t,s,x}^2} & \leq \|f e^{i\Gamma(t,0)\Delta}\|_{L_x^2\to L_{t,x}^2} \|e^{i\Gamma(0,s)\Delta}[f(s)\1_{E_{k'}(s)} F(s)]\|_{L_{s,x}^2} \\
& \leq \|f e^{i\Gamma(t,0)\Delta}\|_{L_x^2\to L_{t,x}^2}\| f \1_{E_{k'}}\|_{L_{s,x}^\infty} \|\1_{E_{k'}}F\|_{L_{s,x}^2} \\
& \leq 2^{\frac{k'}{2}}\|f e^{i\Gamma(t,0)\Delta}\|_{L_x^2\to L_{t,x}^2}\|\1_{E_{k'}} F\|_{L_{t,x}^2}.
\end{align*}
Thus
\begin{align*}
\biggl\| & \int_{\mathcal{R}_{\text{near}}} H_{k,k'}(t,s)\,ds\biggr\|_{L_{t,x}^2} \\ & \leq 2^{\frac{k'}{2}} \|f e^{i\Gamma(t,0)\Delta}\|_{L_x^2\to L_{t,x}^2}  \|\1_{\mathcal{R}_{\text{near}}(t;k,k')}(s)\|_{L_t^\infty L_s^2} \|\1_{E_{k'}} F\|_{L_{t,x}^2}.
\end{align*}

For $k<k'$, we instead estimate
\begin{align*}
\biggl\|& \int_{\mathcal{R}_{\text{near}}} H_{k,k'}(t,s)\,ds\biggr\|_{L_{t,x}^2}\\
&  \leq \|f\1_{E_k}\|_{L_{t,x}^\infty} \biggl\| \int e^{i\Gamma(0,s)\Delta}[f(s) \1_{\mathcal{R}_{\text{near}(t;k,k')}}(s) \1_{E_{k'}(s)}F(s)]\,ds \biggr\|_{L_{t,x}^2} \\
& \leq 2^{\frac{k}{2}}\|[f e^{i\Gamma(t,0)\Delta}]^*\|_{L_{t,x}^2\to L_x^2} \| \1_{\mathcal{R}_{\text{near}}(t;k,k')}(s) \1_{E_{k'}(s)} F(s)\|_{L_{t,s,x}^2} \\
& \leq 2^{\frac{k}{2}} \|f e^{i\Gamma(t,0)\Delta}\|_{L_x^2\to L_{t,x}^2} \| \1_{\mathcal{R}_{\text{near}}(t;k,k')}(s)\|_{L_s^\infty L_t^2} \|\1_{E_{k'}} F\|_{L_{s,x}^2},
\end{align*}
and hence 
\begin{equation}\label{near-temp}
\begin{aligned}
\biggl\| & \int_{\mathcal{R}_{\text{near}}} H_{k,k'}(t,s)\,ds\biggr\|_{L_{t,x}^2}  \\
& \leq 2^{\frac12 [k\wedge k']} \|f e^{i\Gamma(t,0)\Delta}\|_{L_x^2\to L_{t,x}^2} \|\1_{\mathcal{R}_{\text{near}}(t;k,k')}(s)\|_{L_t^\infty L_s^2\cap L_s^\infty L_t^2} \|\1_{E_{k'}}F\|_{L_{t,x}^2},
\end{aligned}
\end{equation}
where $k\wedge k':=\min\{k,k'\}$.  Before estimating the $\1_{\mathcal{R}_{\text{near}}}$ term or addressing the summation in $k,k'$, let us first consider the far region. 

In the far region we may apply the standard dispersive estimate
\[
\| e^{i\Gamma(t,s)\Delta}\|_{L_x^1\to L_x^\infty} \lesssim |\Gamma(t,s)|^{-\frac{d}{2}},\quad s\in\mathcal{R}_{\text{far}}(t;k,k').
\]
Thus
\begin{align*}
\biggl\|& \int_{\mathcal{R}_{\text{far}}} H_{k,k'}(t,s)\,ds\biggr\|_{L_{t,x}^2} \\
& \lesssim\|f\1_{E_k}\|_{L_t^\infty L_x^2} \biggl\| \int |\Gamma(t,s)|^{-\frac{d}{2}} \1_{\mathcal{R}_{\text{far}}(t;k,k')}(s) \|f \1_{E_{k'}} F\|_{L_x^1}\,ds\biggr\|_{L_t^2} \\
& \lesssim 2^{\frac{k}{2}}\|\1_{E_k} \|_{L_t^\infty L_x^2} \biggl\| \int |\Gamma(t,s)|^{-\frac{d}{2}}\1_{\mathcal{R}_{\text{far}}(t;k,k')}(s) \|f \1_{E_{k'}}\|_{L_x^2} \|1_{E_{k'}}(s)F(s)\|_{L_x^2} \,ds\biggr\|_{L_t^2} \\
&\lesssim2^{\frac{k}{2}}2^{\frac{k'}{2}}\|\1_{E_k}\|_{L_t^\infty L_x^2} \|\1_{E_{k'}}\|_{L_t^\infty L_x^2}\biggl\| \int \frac{\1_{\mathcal{R}_{\text{far}}(t;k,k')}(s)}{|\Gamma(t,s)|^{\frac{d}{2}}} \|\1_{E_{k'}}(s)F(s)\|_{L_x^2}\,ds\biggr\|_{L_t^2}.
\end{align*}
Now, observing that
\[
\|\1_{E_k}\|_{L_t^\infty L_x^2} \leq 2^{-\frac{dk}{4}}\|f\|_{L_t^\infty L_x^{d,\infty}}^{\frac{d}{2}}
\]
and defining
\[
\mathcal T_{\text{far};k,k'}\psi(t) = \int |\Gamma(t,s)|^{-\frac{d}{2}}\1_{\mathcal{R}_{\text{far}}(t;k,k')}(s) \psi(s)\,ds,
\]
we find
\begin{align*}
\biggl\|& \int_{\mathcal{R}_{\text{far}}} H_{k,k'}(t,s)\,ds\biggr\|_{L_{t,x}^2}\lesssim2^{-[\frac{d-2}{4}](k+k')}\|f\|_{L_t^\infty L_x^{d,\infty}}^d\|\mathcal{T}_{\text{far};k,k'}\|_{L_t^2\to L_t^2} \|\1_{E_{k'}}F\|_{L_{t,x}^2}.
\end{align*}

In particular, applying Schur's Test, 
\begin{equation}\label{far-temp}
\begin{aligned}
&\biggl\|  \int_{\mathcal{R}_{\text{far}}} H_{k,k'}(t,s)\,ds\biggr\|_{L_{t,x}^2} \\
&\lesssim2^{-[\frac{d-2}{4}](k+k')}\|f\|_{L_t^\infty L_x^{d,\infty}}^d\biggl\|\frac{\1_{\mathcal{R}_{\text{far}}(t;k,k')}(s)}{|\Gamma(t,s)|^{\frac{d}{2}}}\biggr\|_{L_t^\infty L_s^1 \cap L_s^\infty L_t^1}\|\1_{E_{k'}}F\|_{L_{t,x}^2}.
\end{aligned}
\end{equation}

We turn now to the estimation of the quantities involving $\1_{\mathcal{R}_{\text{near}}}$ and $\1_{\mathcal{R}_{\text{far}}}$ appearing in \eqref{near-temp} and \eqref{far-temp}.  First, observe that if $\gamma(t)\equiv 1$, then $\Gamma(t,s)=t-s$, and the problem reduces to the one considered in \cite{Hadama} (with $U(t)=e^{it\Delta}$, $H=L^2$, and $\sigma=\tfrac{d}{2}$ in the notation of \cite{Hadama}). In this case, the proof in \cite{Hadama} makes use of the explicit bounds
\begin{equation}\label{gamma=1case}
\begin{aligned}
& \|\1_{\mathcal{R}_{\text{near}}(t;k,k')}(s)\|_{L_t^\infty L_s^2 \cap L_s^\infty L_t^2}  \sim |\ell(k,k')|^{\frac12}\quad (\text{if }\gamma(t)\equiv 1),\\
& \|\1_{\mathcal{R}_{\text{far}}(t;k,k')}(s)|\Gamma(t,s)|^{-\frac{d}{2}}\|_{L_t^\infty L_s^1 \cap L_s^\infty L_t^1} \sim |\ell(k,k')|^{1-\frac{d}{2}}\quad (\text{if }\gamma(t)\equiv 1). 
\end{aligned}
\end{equation}
We will prove that for admissible dispersion maps, one can obtain \eqref{gamma=1case} up to some multiplicative constant:

\begin{proposition}\label{P} There exists $C=C(\gamma)>0$ such that
\begin{align*}
\|\1_{\mathcal{R}_{\text{near}}(t;k,k')}(s)\|_{L_t^\infty L_s^2 \cap L_s^\infty L_t^2}  &\leq  C|\ell(k,k')|^{\frac12}, \\
 \|\1_{\mathcal{R}_{\text{far}}(t;k,k')}(s)|\Gamma(t,s)|^{-\frac{d}{2}}\|_{L_t^\infty L_s^1 \cap L_s^\infty L_t^1} &\leq C |\ell(k,k')|^{1-\frac{d}{2}}.
\end{align*}
\end{proposition}

\begin{proof}  Both estimates will be derived from a generalization of the argument used in \cite[Proof of Theorem~2]{MVH}, which derived Strichartz estimates for $e^{i\Gamma(t,s)\Delta}$ from the standard Strichartz estimates for $e^{it\Delta}$. 

We begin with the near estimate.  Using Cauchy--Schwarz, it suffices to show the following: 
\[
\exists C>0:\forall\ell>0\ \forall t\in\R \quad |\{s\in\R: s<t\qtq{and} |\Gamma(t,s)|\leq\ell\}| \leq C\ell. 
\]
To this end, we fix $\ell>0$, $t\in\R$, and first write
 \[
|\Gamma(t,s)|=|\Gamma(s)-\Gamma(t)|,\qtq{where}\Gamma(\cdot)\qtq{denotes}\Gamma(\cdot,0).
\] 
We begin by writing
\[
|\{s:s<t\qtq{and} |\Gamma(t,s)|\leq\ell\}| = \int_{-\infty}^t \chi(\Gamma(s))\,ds \leq \int_\R \chi(\Gamma(s))\,ds,
\]
where
\[
\chi(\sigma) = \begin{cases} 1 & \text{if } |\sigma - \Gamma(t)| \leq \ell, \\ 0 & \text{otherwise}.\end{cases}
\]
Now (following \cite{MVH}) split $\R$ into disjoint sets $\mathbb{I}$ and $\mathbb{D}$ on which $\Gamma(\cdot)$ is increasing and decreasing, respectively.  Each of these sets is a disjoint union of intervals, and if $N_\gamma$ denotes the number of discontinuities of $\gamma$ in one period, then the number of intervals comprising $\mathbb{I}$ or $\mathbb{D}$ within any unit interval is at most $N_\gamma+1$. 

It suffices to prove a bound on the integral over $\mathbb{I}$.  To this end, for $n\in\mathbb{Z}$ we decompose
\[
\Gamma^{-1}([n,n+1))\cap\mathbb{I} = \bigcup_{k=1}^{K_n} U_k^n,
\]
where $\{U_k^n\}_{k=1}^{K_n}$ is disjoint collection of intervals, with each $U_k^n$ contained in a distinct subinterval of $\mathbb{I}$. Arguing as in \cite{MVH}, the key observation is that admissibility of $\gamma$ guarantees the bound
\[
\sup_{n\in\mathbb{Z}} K_n \leq K_\gamma := 1+N_\gamma\bigl\{ 1+\tfrac{1+4\|\gamma\|_{L^\infty}}{|\langle\gamma\rangle|}\bigr\}
\]
(see \cite[(2.6)]{MVH}).  As $s\mapsto \Gamma(s)$ is injective on each interval $U_k^n$, we can therefore make the change of variables $\sigma=\Gamma(s)$ on each $U_k^n$ and obtain the following estimate: 
\begin{align*}
\int_\mathbb{I} \chi(\Gamma(s))\,ds &  \leq \sum_{n\in\mathbb{Z}} \sum_{k=1}^{K_n} \int_{U_k^n} \chi(\Gamma(s))\,ds \\
& \leq K_\gamma \|\gamma^{-1}\|_{L^\infty} \sum_{n\in\mathbb{Z}} \int_n^{n+1} \chi(\sigma)\,d\sigma  \\
& \leq 2K_\gamma \|\gamma^{-1}\|_{L^\infty}\, \ell,
\end{align*}
which is acceptable. As a similar argument treats the contribution of $\mathbb{D}$, this completes the proof of the near estimate.  

We turn to the far estimate. It will suffice to prove that there exists $C>0$ such that for any $t\in\R$ and any $\ell>0$,  
\[
\int_\R \Phi(\Gamma(s))\,ds \leq C|\ell|^{1-\frac{d}{2}},
\]
where
\[
\Phi(\sigma) = \begin{cases} |\sigma-\Gamma(t)|^{-\frac{d}{2}} & \text{if }|\sigma-\Gamma(t)|>\ell, \\ 0 & \text{otherwise}.\end{cases}
\]
We estimate this integral using the same strategy as above, decomposing $\R=\mathbb{I}\cup\mathbb{D}$ and $\Gamma^{-1}([n,n+1))\cap\mathbb{I}=\bigcup_{k=1}^{K_n} U_k^n$.  Then
\begin{align*}
\int_{\mathbb{I}} \Phi(\Gamma(s))\,ds &\leq \sum_{n\in\mathbb{Z}}\sum_{k=1}^{K_n} \int_{U_n^k} \Phi(\Gamma(s))\,ds \\
& \leq K_\gamma\|\gamma^{-1}\|_{L^\infty} \sum_{n\in\mathbb{Z}} \int_n^{n+1} \Phi(\sigma)\,d\sigma \\
& \leq K_\gamma\|\gamma^{-1}\|_{L^\infty} \int_{|\sigma|>\ell} |\sigma|^{-\frac{d}{2}}\,d\sigma \\
& \leq \tfrac{4}{d-2} K_\gamma\|\gamma^{-1}\|_{L^\infty} |\ell|^{1-\frac{d}{2}},
\end{align*}
which is acceptable.  As a similar argument treats the contribution of $\mathbb{D}$, this completes the proof of the far estimate. \end{proof}

Having established the analogues of \eqref{gamma=1case} for the case of an admissible dispersion map, it is clear in principle that the rest of the proof in \cite{Hadama} will go through.  We turn to the details and return to the proof of \eqref{de1}.

Combining \eqref{l1l2} with the near and far estimates \eqref{near-temp} and \eqref{far-temp}, and inserting the bounds established in Proposition~\ref{P}, we obtain
\begin{align*}
\biggl\| & \int_{-\infty}^t f(t) e^{i\Gamma(t,s)\Delta}[f(s)F(s)]\,ds \biggr\|_{L_{t,x}^2} \\
& \lesssim \|fe^{i\Gamma(t,0)\Delta}\|_{L_x^2\to L_{t,x}^2} \bigl\| 2^{\frac12[k\wedge k']} [\ell(k,k')]^{\frac12} \1_{E_k'} F\bigr\|_{\ell_k^2 \ell_{k'}^1 L_{t,x}^2} \\
& \quad + \|f\|_{L_t^\infty L_x^{d,\infty}}^d \bigl\|2^{-[\frac{d-2}{4}](k+k')} [\ell(k,k')]^{1-\frac{d}{2}} \1_{E_{k'}}F\bigr\|_{\ell_k^2 \ell_{k'}^1 L_{t,x}^2}. 
\end{align*}
We now optimize in the choice of $\ell(k,k')$, choosing
\begin{align*}
\ell(k,k') =\bigl[\|f e^{i\Gamma(t,0)\Delta}\|_{L_x^2\to L_{t,x}^2}^{-\frac{2}{d-1}} \|f\|_{L_t^\infty L_x^{d,\infty}}^{\frac{2d}{d-1}} \bigr]\cdot 2^{-\frac{1}{d-1}[k\wedge k']} 2^{-\frac{d-2}{2(d-1)}[k+k']}. \end{align*}
Writing
\[
M=\|f e^{i\Gamma(t,0)\Delta}\|_{L_x^2\to L_{t,x}^2}^{\frac{d-2}{d-1}} \|f\|_{L_t^\infty L_x^{d,\infty}}^{\frac{d}{d-1}}
\]
and using Young's convolution inequality, we therefore obtain
\begin{align*}
\biggl\|  \int_{-\infty}^t & f(t) e^{i\Gamma(t,s)\Delta}[f(s)F(s)]\,ds \biggr\|_{L_{t,x}^2} \\
& \lesssim M\| 2^{\frac{d-2}{2(d-1)}[k\wedge k']}2^{-\frac{d-2}{4(d-1)}[k+k']} \1_{E_{k'}} F\|_{\ell_k^2\ell_{k'}^1 L_{t,x}^2} \\
& \lesssim M \| 2^{-\frac{d-2}{4(d-1)}|k-k'|} \1_{E_{k'}} F\|_{\ell_k^2 \ell_{k'}^1 L_{t,x}^2} \\
& \lesssim M \| 2^{-\frac{d-2}{4(d-1)}|k|} \ast \| \1_{E_k} F\|_{L_{t,x}^2}\|_{\ell_k^2} \\
& \lesssim M \| \1_{E_k} F\|_{\ell_k^2 L_{t,x}^2}.
\end{align*}
Now, by the disjointness of the sets $E_k$,
\begin{align*}
\|\1_{E_k} F\|_{\ell_k^2 L_{t,x}^2}^2 \leq \iint \biggl(\sum_k \1_{E_k}(t,x)\biggr)|F(t,x)|^2\,dx\,dt \leq \|F\|_{L_{t,x}^2}^2.
\end{align*}
Thus, applying Young's inequality,
\begin{align*}
\biggl\|  \int_{-\infty}^t & f(t) e^{i\Gamma(t,s)\Delta}[f(s)F(s)]\,ds \biggr\|_{L_{t,x}^2} \\
& \lesssim \|f e^{i\Gamma(t,0)\Delta}\|_{L_x^2\to L_{t,x}^2}^{\frac{d-2}{d-1}} \|f\|_{L_t^\infty L_x^{d,\infty}}^{\frac{d}{d-1}} \|F\|_{L_{t,x}^2} \\
& \leq \biggl[\tfrac14\|f e^{i\Gamma(t,0)\Delta}\|_{L_x^2 \to L_{t,x}^2}^2 + C\|f\|_{L_t^\infty L_x^{d,\infty}}^2 \biggr]\|F\|_{L_{t,x}^2}. 
\end{align*}

Writing
\[
A := \sup\biggl\{\biggl\| \int_{-\infty}^t f(t) e^{i\Gamma(t,s)\Delta}[f(s)F(s)]\,ds \biggr\|_{L_{t,x}^2}: \|F\|_{L_{t,x}^2} = 1\biggr\},
\]
we therefore obtain
\begin{equation}\label{almost-there}
A \leq \tfrac14 \|f e^{i\Gamma(t,0)\Delta}\|_{L_x^2\to L_{t,x}^2}^2 + C\|f\|_{L_t^\infty L_x^{d,\infty}}^2. 
\end{equation}
The final observation is that we can write
\begin{align*}
[fe^{i\Gamma(t,0)\Delta}][fe^{i\Gamma(t,0)\Delta}]^* = P + P^*,
\end{align*}
where $P$ is the principal operator under consideration, i.e. 
\[
[PF](t) = \int_{-\infty}^t f(t) e^{i\Gamma(t,s)\Delta}[f(s)F(s)]\,ds.
\]
Thus, by the $TT^*$ estimate, we can continue from \eqref{almost-there} to obtain
\[
A \leq \tfrac12 A + C\|f\|_{L_t^\infty L_x^{d,\infty}}^2.
\]
Noting that \eqref{form} guarantees $A<\infty$, we therefore obtain  
\[
A\leq 2C\|f\|_{L_t^\infty L_x^{d,\infty}}^2,
\]
which completes the proof of \eqref{de1} in the reduced setting of $f=g\geq 0$ satisfying \eqref{form}.  As described above, the rest of the argument in \cite{Hadama} goes through unchanged, and thus we complete the proof of Theorem~\ref{T}.

\end{document}